\documentclass[11pt,reqno]{amsart}

\usepackage{amsmath,amsthm,amssymb,mathtools}
\usepackage{booktabs,longtable,array}
\usepackage[margin=1in]{geometry}
\usepackage{microtype}
\usepackage[colorlinks=true,linkcolor=blue,citecolor=blue,urlcolor=blue]{hyperref}
\usepackage{orcidlink}

\allowdisplaybreaks
\newtheorem{theorem}{Theorem}[section]
\newtheorem{proposition}[theorem]{Proposition}
\newtheorem{lemma}[theorem]{Lemma}
\newtheorem{corollary}[theorem]{Corollary}
\theoremstyle{definition}

\newtheorem{example}[theorem]{Example}
\theoremstyle{remark}
\newtheorem{remark}[theorem]{Remark}

\newcommand{\A}{\mathbb A}
\newcommand{\C}{\mathbb C}

\newcommand{\Q}{\mathbb Q}
\newcommand{\Z}{\mathbb Z}
\newcommand{\GL}{\operatorname{GL}}
\newcommand{\SL}{\operatorname{SL}}
\newcommand{\Sym}{\operatorname{Sym}}
\newcommand{\Ind}{\operatorname{Ind}}
\newcommand{\AI}{\operatorname{AI}}

\newcommand{\spc}{\operatorname{sp}}

\newcommand{\ch}{\operatorname{ch}}

\newcommand{\GammaR}{\Gamma_{\mathbb R}}
\newcommand{\GammaC}{\Gamma_{\mathbb C}}
\newcommand{\eps}{\varepsilon}
\newcommand{\qbinom}[2]{\genfrac{[}{]}{0pt}{}{#1}{#2}}

\title{An evident corollary arising from Newton--Thorne}

\author{Shenghao Hua~\orcidlink{0000-0002-7210-2650}}
\address[1]{Shanghai Institute for Mathematics and Interdisciplinary Sciences (SIMIS), Shanghai, 200433, China}
\address[2]{Research Institute of Intelligent Complex Systems, Fudan University, Shanghai, 200433, China}
\email{shenghao.hua@simis.cn}


\hypersetup{
 pdftitle={An evident corollary arising from Newton--Thorne},
 pdfauthor={Shenghao Hua},
 pdfkeywords={elliptic curve, symmetric power, tensor product, isobaric sum,
 Weil--Deligne representation}
}

\begin{document}

\begin{abstract}
Let $f$ be a primitive holomorphic newform of weight at least $2$ and
arbitrary level.  We consider the representations obtained from its
symmetric powers by taking tensor products, symmetric powers, and isobaric
sums.  Classical $\GL_2$ representation theory decomposes each such
expression into symmetric powers with determinant twists.  This gives an
automorphic realization and a factorization of the associated $L$-function.
At a prime dividing the level, the same decomposition must be applied to the
full Weil--Deligne parameter before inertia invariants and the kernel of
monodromy are taken.  For elliptic curves over $\Q$, we record the resulting
local factors, conductors, and root numbers using the calculations of
Dummigan--Martin--Watkins and Martin--Watkins.
\end{abstract}

\maketitle
\begingroup
\small
\tableofcontents
\endgroup

\section{Introduction}

Langlands functoriality predicts transfers of automorphic representations associated with homomorphisms of $L$-groups.
For $\GL_2$, symmetric-power lifts are among the best understood examples.
Gelbart and Jacquet constructed the symmetric-square lift, Kim and Shahidi
the symmetric-cube lift, and Kim the symmetric-fourth lift \cite{GelbartJacquet,KimShahidi,Kim}.  Newton and Thorne subsequently proved the automorphy, and indeed the cuspidality, of every positive symmetric power of a non-CM holomorphic newform of weight at least two and arbitrary level \cite{NewtonThorne2021a,NewtonThorne2021b}.  For a CM newform, the symmetric powers decompose into Hecke characters of $\Q$ and automorphic inductions of Hecke characters from the CM field.

The general $\GL_m\times\GL_n$ theory of Jacquet, Piatetski-Shapiro, and Shalika develops the local factors and global integral representations for the Rankin--Selberg $L$-functions of cuspidal automorphic representations \cite{JPSS}.
This is distinct from a general automorphic transfer $\pi\boxtimes\pi'$ to $\GL_{mn}$, which is not known in arbitrary ranks.

Let $f$ be a primitive holomorphic cusp form of weight $k\geq2$, level $N$, and nebentypus $\chi$, and let $\pi_f^u$ be its unitary automorphic representation.
Let $V$ be the standard representation of $\GL_2(\C)$.
We consider polynomial representations obtained from the representations $\Sym^nV$ by direct sums, tensor products, and symmetric powers.
Every irreducible constituent is of the form
\begin{equation}\label{eq:intro-atoms}
 R_{a,b}=\Sym^aV\otimes\det(V)^b,
 \qquad a,b\geq0.
\end{equation}
Evaluating this decomposition on the parameter of $f$ gives an isobaric sum of twists of the known symmetric-power lifts.
In particular, no general higher-rank tensor-product or symmetric-power transfer is needed.

At an unramified prime, the local factors can be computed directly from the two Satake parameters.
This is no longer true at a ramified prime, since the standard Euler factor retains only the action of Frobenius on $(\ker N)^{I_p}$.

For example, suppose that an elliptic curve $E/\Q$ has additive potentially multiplicative reduction at $p$.
Then $L_p(s,E)=1$, whereas
\begin{equation}\label{eq:intro-counterexample}
 L_p(s,\Sym^2 E)=(1-p^{-s})^{-1}.
\end{equation}
Thus one must first apply the algebraic operation to the full Weil--Deligne
representation, including its monodromy operator, and only afterwards form the local Euler factor.

Call $F$ admissible if it is obtained from the representations $\Sym^nV$, $n\geq0$, by finitely many direct sums, tensor products, and operations $\Sym^m$, $m\geq0$, with $\Sym^0W=1$.
We write $\Lambda$ for the product of the archimedean factor and all finite local factors, omitting the conductor power.
The normalization is given in Section~\ref{sec:normalizations}.

\begin{theorem}\label{thm:global}
Let $f$ be a primitive holomorphic cusp form of weight $k\geq2$, arbitrary level, and nebentypus $\chi$.
For every nonzero admissible expression $F$,
there is a unique finitely supported array
$c_F(a,b)\in\Z_{\geq0}$ such that
\begin{equation}\label{eq:global-algebra-decomp}
 F(V)\simeq
 \bigoplus_{a,b\geq0}R_{a,b}^{\oplus c_F(a,b)}.
\end{equation}
Its unitary automorphic realization is
\begin{equation}\label{eq:unitary-lift}
 \Pi_F^u=
 \boxplus_{a,b\geq0}
 \left(\Sym^a\pi_f^u\otimes\chi^b\right)^{\boxplus c_F(a,b)}.
\end{equation}
At every finite prime, including every prime dividing the level,
\begin{equation}\label{eq:all-local-product}
 L_p(s,F(f))=
 \prod_{a,b\geq0}
 L_p\left(s-b(k-1),\Sym^af\otimes\chi^b\right)^{c_F(a,b)}.
\end{equation}
The same identity holds for the archimedean-completed factors
\begin{equation}\label{eq:completed-product}
 \Lambda(s,F(f))=
 \prod_{a,b\geq0}
\Lambda\left(s-b(k-1),\Sym^af\otimes\chi^b\right)^{c_F(a,b)}.
\end{equation}
If $F$ has constituents of different polynomial degrees, this is a constituentwise motivic product.
It need not have a single motivic normalization or a single functional-equation center.

At a ramified prime, formula \eqref{eq:all-local-product} means that $F$ is first evaluated on the complete Weil--Deligne parameter and only then is
$(\ker N)^{I_p}$ taken.
It is not obtained by applying $F$ to the roots visible in $L_p(s,f)$.

For an elliptic curve $E/\Q$, let $\tau_{E,p}$ be its arithmetic Weil--Deligne parameter and let $P_{p,F}(T)$ be the Euler polynomial of $F(\tau_{E,p})$.
Write
\begin{equation*}
 P_{p,a}(T)=\det\left(1-r_{p,a}(\Phi_p)T
 \mathrel{\big|}(\ker N_{p,a})^{I_p}\right),
 \qquad (r_{p,a},N_{p,a})=\Sym^a\tau_{E,p}.
\end{equation*}
Then
\begin{equation}\label{eq:intro-local-assembly}
 P_{p,F}(T)=
 \prod_{a,b\geq0}P_{p,a}(p^bT)^{c_F(a,b)}.
\end{equation}
If $f_{p,a}$ and $W_{p,a}$ denote the conductor exponent and local root number of $\Sym^aE$, then
\begin{equation}\label{eq:intro-conductor-root}
 f_p(F(E))=\sum_{a,b}c_F(a,b)f_{p,a},
 \qquad
 W_p(F(E))=\prod_{a,b}W_{p,a}^{c_F(a,b)}.
\end{equation}
\end{theorem}

Theorem~\ref{thm:global} combines the symmetric-power theorem of Newton and
Thorne, and automorphic induction in the CM case, with the classical $\GL_2$ decomposition into the representations $R_{a,b}$.
The identities of local factors then follow from the additivity of Weil--Deligne
representations and the multiplicativity of local factors.

Section~\ref{sec:normalizations} fixes the two normalizations and proves Theorem~\ref{thm:global}.  Section~\ref{sec:ramified} treats the ramified Weil--Deligne calculation.  Appendix~\ref{app:algebra} gives the
Clebsch--Gordan and Cayley--Sylvester multiplicities and the resulting recursive classification.
Appendix~\ref{app:local} specializes to an
elliptic curve over $\Q$ and records the local factors, conductors, and root numbers obtained from \cite{DMW,MartinWatkins}, including the wild cases at $2$ and $3$.
The final appendix treats CM curves and the archimedean place.

\section{The global construction}

\subsection{Two normalizations}\label{sec:normalizations}

Let $\pi_f^u$ denote the unitary automorphic representation of $\GL_2(\A_\Q)$ attached to $f$.
Its central character is the finite-order nebentypus $\chi$.
Thus the determinant of its local Langlands parameter $\phi_{f,v}$ is $\chi_v$.

For arithmetic Euler factors, let
\begin{equation}\label{eq:motivic-WD}
 \tau_{f,p}=(r_{f,p},N_{f,p}),\qquad
 \det r_{f,p}=\chi_p\omega_p^{-(k-1)},\qquad
 \omega_p(\Phi_p)=p^{-1},
\end{equation}
where $\Phi_p$ is geometric Frobenius.
At an unramified prime the eigenvalues $\alpha_p,\beta_p$ of $r_{f,p}(\Phi_p)$ satisfy
\begin{equation*}
 \alpha_p+\beta_p=a_p(f),\qquad
 \alpha_p\beta_p=\chi(p)p^{k-1}.
\end{equation*}
The relation between the two normalizations is
\begin{equation}\label{eq:normalization-relation}
 \phi_{f,p}=\tau_{f,p}\otimes\omega_p^{(k-1)/2}.
\end{equation}
The normalizations may be compared at a good prime as follows.
\begin{center}
\renewcommand{\arraystretch}{1.15}
\begin{tabular}{@{}llll@{}}
\toprule
Parameter & Determinant & Product of roots & Role \\
\midrule
$\phi_{f,p}$ & $\chi_p$ & $\chi(p)$ & unitary isobaric decomposition \\
$\tau_{f,p}$ & $\chi_p\omega_p^{-(k-1)}$
& $\chi(p)p^{k-1}$ & classical Euler factors \\
\bottomrule
\end{tabular}
\end{center}

Let $V$ be the standard representation of $\GL_2(\C)$.
For $a,b\geq0$ put
\begin{equation}\label{eq:Rab}
 R_{a,b}(V)=\Sym^aV\otimes\det(V)^b.
\end{equation}
It has dimension $a+1$ and polynomial degree $a+2b$.
On the unitary parameter,
\begin{equation}\label{eq:Rab-unitary-general}
 R_{a,b}(\phi_{f,v})
 =\Sym^a\phi_{f,v}\otimes\chi_v^b.
\end{equation}
On the motivic parameter,
\begin{equation}\label{eq:Rab-motivic}
 R_{a,b}(\tau_{f,p})
 =\Sym^a\tau_{f,p}\otimes\chi_p^b\omega_p^{-b(k-1)}.
\end{equation}
For an elliptic curve one has $k=2$ and $\chi=1$.
In that case the determinant twist disappears in the unitary normalization and becomes the Tate shift $s\mapsto s-b$ in the arithmetic normalization.

We write $\Lambda(s,\Sym^af\otimes\chi^b)$ for the product of all finite local factors and the archimedean factor in this arithmetic (motivic) normalization.
Thus $\Lambda$ here is archimedean-completed and includes the bad-prime factors, but it does not include the conventional conductor power $N^{s/2}$.
This convention makes shifted products literal, with no irrelevant nonzero conductor constants.

\subsection{Admissible expressions}

An admissible expression is obtained from the representations $\Sym^nV$,
$n\geq0$, by a finite number of the following operations.
\begin{enumerate}
 \item A direct sum.
 \item A tensor product.
 \item A symmetric power $\Sym^m$, with $m\geq0$.
\end{enumerate}
We use $\Sym^0(W)=1$.  For algebraic bookkeeping we also adjoin the zero
representation as the empty direct sum, but the zero object is not called an
automorphic representation.  Since $\Sym^0V=1$ is already an initial object,
allowing the operation $\Sym^0$ is redundant but convenient in formulas.
We distinguish throughout between irreducible constituents that can occur
and complete representations produced by the allowed operations.
Proposition~\ref{prop:all-atoms} classifies the former;
Proposition~\ref{prop:whole-expressions} classifies the latter.
When the expression is evaluated automorphically, a direct sum of parameters
is represented by the isobaric sum, or equivalently by the Langlands quotient
of the corresponding normalized parabolic induction.  We do not identify an
arbitrary reducible induced representation with this quotient.

\subsection{Proof and consequences of the main theorem}

\begin{proof}[Proof of Theorem~\ref{thm:global}]
Every irreducible polynomial representation of $\GL_2(\C)$ has the form
$R_{a,b}$ in \eqref{eq:Rab}.  Complete reducibility gives
\eqref{eq:global-algebra-decomp}, and highest weights give uniqueness.  The
standard highest-weight classification used here may be found in
\cite[Sections~6.1--6.2]{FultonHarris}.  The formulae in
Appendix~\ref{app:algebra} calculate every multiplicity.

For non-CM $f$, Theorem~A of Newton--Thorne gives the cuspidal automorphic
representation $\Sym^a\pi_f^u$ for every $a\geq1$
\cite[Theorem~A]{NewtonThorne2021b}.  Local--global compatibility for the
attached Galois representations identifies its finite-place
Frobenius-semisimplified Weil--Deligne parameters with the symmetric powers of those of
$f$, including at supercuspidal bad places.  In the elliptic-curve case this
is compatible with the completed symmetric-power function of Dummigan--Martin--Watkins used in
\cite[Corollary~B]{NewtonThorne2021b}.  For CM $f$,
the same conclusion follows from automorphic induction.  In weight $2$ the
explicit decomposition is given in Proposition~\ref{prop:CM}.  Twisting by
$\chi^b$ and taking an isobaric sum preserve automorphy.  This constructs
\eqref{eq:unitary-lift}.  Evaluating
\eqref{eq:global-algebra-decomp} on $\phi_{f,v}$ and using
\eqref{eq:Rab-unitary-general} gives its local parameter.

For the motivic normalization, evaluate the same algebraic identity on the
complete Weil--Deligne parameter $\tau_{f,p}$.  Equation
\eqref{eq:Rab-motivic} gives the finite-order twist $\chi_p^b$ and replaces
$s$ by $s-b(k-1)$.  Local $L$-factors are multiplicative under direct sums,
which proves \eqref{eq:all-local-product}.  Local compatibility at infinity
proves \eqref{eq:completed-product}.

For $E/\Q$ one has $\chi=1$ and
$\det\tau_{E,p}=\omega_p^{-1}$.  Hence the summand $R_{a,b}$ replaces
$P_{p,a}(T)$ by $P_{p,a}(p^bT)$.  Artin conductors are additive and epsilon
factors are multiplicative under direct sums.  The unramified positive real
twist $\omega_p^{-b}$ changes neither the conductor exponent nor the phase of
the local root number.  This proves \eqref{eq:intro-local-assembly} and
\eqref{eq:intro-conductor-root}.
\end{proof}

\begin{corollary}[Purity and cuspidality]\label{cor:purity-cuspidality}
If $F$ is homogeneous of polynomial degree $D$, then
\begin{equation}\label{eq:homogeneous-shift}
 F(\tau_f)=F(\phi_f)\otimes\omega^{-D(k-1)/2},
 \qquad
 L(s,F(f))=L\left(s-\frac{D(k-1)}2,\Pi_F^u\right).
\end{equation}
If $F$ contains several polynomial degrees, equations
\eqref{eq:all-local-product} and \eqref{eq:completed-product} are instead
blockwise motivic products.  There is then no single shift relating the
whole product to the standard $L$-function of $\Pi_F^u$.

If $f$ is non-CM, every $\Sym^a\pi_f^u$ with $a\geq1$ is cuspidal.  Hence
an output of total rank at least two is cuspidal precisely when, after equal
constituents in \eqref{eq:unitary-lift} are combined, there is one
multiplicity-one block with $a\geq1$.  A rank-one output is a single Hecke
character and is cuspidal under the standard convention for $\GL_1$.
\end{corollary}

\begin{proof}
The first assertion follows by applying a homogeneous polynomial functor of
degree $D$ to \eqref{eq:normalization-relation}.  The cuspidality assertion
follows from Newton--Thorne in the non-CM case and from the definition of an
isobaric sum.
\end{proof}

\begin{remark}\label{rem:no-general-functoriality}
The notation $F(\pi_f)$ does not assert a general functorial tensor product or
a general symmetric-power lift for a representation of $\GL_n$.  It denotes
the isobaric representation obtained after the two-dimensional parameter has
been decomposed as in \eqref{eq:global-algebra-decomp}.
\end{remark}

\begin{corollary}[Elliptic-curve outputs]\label{cor:all-unitary}
Let $E/\Q$ be an elliptic curve.  As $F$ ranges over all admissible
expressions, its unitary automorphic outputs
are exactly
\begin{equation}\label{eq:all-unitary}
 \boxplus_{a\geq0}
 \left(\Sym^a\pi_E^u\right)^{\boxplus C(a)},
 \qquad C(a)\in\Z_{\geq0},
\end{equation}
where $C(a)=0$ for all but finitely many $a$, and the array is not identically
zero.
\end{corollary}

\begin{proof}
Theorem~\ref{thm:global} shows that every output has this form.  Conversely,
all $\Sym^a\pi_E^u$ are among the initial objects, and finite isobaric sums are
allowed.  We include $\Sym^0\pi_E^u=1$.
\end{proof}

Here $C(a)=\sum_b c_F(a,b)$.  Passing to the unitary output forgets the
index $b$; consequently, \eqref{eq:all-unitary} does \emph{not} determine
the arithmetic $L$-function.  Its Tate shifts require the full array
$c_F(a,b)$ in \eqref{eq:all-local-product}.

\section{Ramified places and the Weil--Deligne parameter}\label{sec:ramified}

The purpose of this section is to isolate the point at which the unramified
Satake-root calculation stops being valid.  The formal statement applies to
every holomorphic newform.  In the explicit examples and all subsequent local
tables we specialize to an elliptic curve $E/\Q$.  Thus $k=2$, $\chi=1$, and
its motivic parameter $\tau_p$ satisfies \eqref{eq:motivic-WD} with
$\det r_p=\omega_p^{-1}$.

For a Weil--Deligne representation $\tau=(r,N)$ over $\Q_p$, define
\begin{equation}\label{eq:Dfunctor}
 D_p(\tau)=(\ker N)^{I_p}.
\end{equation}
Its Euler polynomial and local factor are
\begin{equation}\label{eq:euler-definition}
 P_{p,\tau}(T)=
 \det\left(1-r(\Phi_p)T\mid D_p(\tau)\right),
 \qquad
 L_p(s,\tau)=P_{p,\tau}(p^{-s})^{-1}.
\end{equation}
These are the standard Weil--Deligne conventions; see, for example,
\cite{Rohrlich}.
We use geometric Frobenius and the standard additive character of $\Q_p$.
The local root number $W_p(\tau)$ is the unit-modulus phase of the associated
epsilon factor after unitary normalization.  For the self-dual symmetric
powers considered below it lies in $\{\pm1\}$.

\begin{proposition}[What survives at a bad prime]\label{prop:WD-decomposition}
Suppose
\begin{equation*}
 F(V)\simeq
 \bigoplus_{a,b}R_{a,b}(V)^{\oplus c(a,b)}.
\end{equation*}
Then, for every Weil--Deligne parameter $\tau_p$ of an elliptic curve,
\begin{equation}\label{eq:WD-F-decomp}
 F(\tau_p)\simeq
 \bigoplus_{a,b}
 \left(\Sym^a\tau_p\otimes(\det\tau_p)^b\right)^{\oplus c(a,b)}.
\end{equation}
However, in general
\begin{equation}\label{eq:D-not-monoidal}
 D_p(\tau_1\otimes\tau_2)\neq D_p(\tau_1)\otimes D_p(\tau_2),
 \qquad
 D_p(\Sym^m\tau)\neq\Sym^mD_p(\tau).
\end{equation}
Consequently, the roots visible in $P_{p,\tau_p}$ cannot be used as ramified
Satake parameters.
\end{proposition}

\begin{proof}
On a tensor product the monodromy operator is
\begin{equation*}
 N_1\otimes1+1\otimes N_2.
\end{equation*}
On a Schur functor it is the image of $N$ under the differential of that
functor.  Every isomorphism used to decompose $F(V)$ is $\GL_2$-equivariant.
It therefore intertwines both the Weil action and the differentiated action
of $N$.  This proves \eqref{eq:WD-F-decomp}.

Taking the kernel of $N$ and then inertia invariants is additive, but it is not
a tensor functor.  The examples below prove both failures in
\eqref{eq:D-not-monoidal} explicitly.
\end{proof}

\begin{example}[Split multiplicative reduction]\label{ex:split-mult}
Put $T=p^{-s}$.  One has
\begin{equation*}
 P_{p,E}(T)=1-T.
\end{equation*}
Nevertheless,
\begin{equation*}
 \tau_p\otimes\tau_p\simeq\Sym^2\tau_p\oplus\det\tau_p
\end{equation*}
and hence
\begin{equation}\label{eq:split-tensor-counter}
 P_{p,E\otimes E}(T)=(1-T)(1-pT).
\end{equation}
The second root comes from the determinant line and cannot be obtained by
squaring the only root visible in $P_{p,E}$.
\end{example}

\begin{example}[Two representations with the same standard factor]
If $E$ has additive potentially multiplicative reduction, then
\begin{equation*}
 P_{p,E}(T)=1,
 \qquad
 P_{p,\Sym^2E}(T)=1-T.
\end{equation*}
If $E$ has potentially good reduction and inertia $C_2$ acts by $-1$ on the
two-dimensional space, then again $P_{p,E}(T)=1$.  In the second case inertia
acts trivially on the symmetric square, so
\begin{equation}\label{eq:C2-counter}
 P_{p,\Sym^2E}(T)=
 (1-\gamma_p^2T)(1-pT)(1-\delta_p^2T).
\end{equation}
Here $\gamma_p$ and $\delta_p$ are the two eigenvalues of Frobenius on the
\emph{full} two-dimensional Weil--Deligne representation; they are not roots
visible in $P_{p,E}$, which has no roots in this example.
Thus even the standard Euler polynomial does not determine the
symmetric-square Euler polynomial.
\end{example}

The same calculation gives a compact comparison for
$E\otimes E=\Sym^2E\oplus\det E$.  In the good row
$\alpha_p\beta_p=p$; in the potentially good row
$\gamma_p\delta_p=p$; and $T=p^{-s}$.
\begin{center}
\renewcommand{\arraystretch}{1.15}
\begin{tabular}{@{}lll@{}}
\toprule
Local type & $P_{p,E}(T)$ & $P_{p,E\otimes E}(T)$ \\
\midrule
Good & $(1-\alpha_pT)(1-\beta_pT)$
& $(1-\alpha_p^2T)(1-pT)^2(1-\beta_p^2T)$ \\
Split multiplicative & $1-T$ & $(1-T)(1-pT)$ \\
Additive potentially multiplicative & $1$ & $(1-T)(1-pT)$ \\
Potentially good, $I=C_2$ & $1$
& $(1-\gamma_p^2T)(1-pT)^2(1-\delta_p^2T)$ \\
\bottomrule
\end{tabular}
\end{center}
The last two rows have the same standard polynomial but different tensor
square polynomials.  This is the simplest concrete reason to retain the full
local parameter.

For later reference, let
\begin{equation}\label{eq:Ppa}
 P_{p,a}(T)=
 \det\left(
 1-r_{p,a}(\Phi_p)T
 \mathrel{\big|}(\ker N_{p,a})^{I_p}
 \right),
 \qquad
 (r_{p,a},N_{p,a})=\Sym^a\tau_p.
\end{equation}
Proposition~\ref{prop:WD-decomposition} gives the correct assembly formula
\begin{equation}\label{eq:P-F-assembly}
 P_{p,F}(T)=
 \prod_{a,b}P_{p,a}(p^bT)^{c_F(a,b)}.
\end{equation}
The polynomials $P_{p,a}$ must be calculated from the complete local type.
They are listed in Appendix~\ref{app:local}.

\begin{proposition}[Conductor and root number assembly]\label{prop:conductor-assembly}
Let
\begin{equation*}
 f_{p,a}=a(\Sym^a\tau_p),\qquad
 W_{p,a}=
 \frac{\eps(\tfrac12,\Sym^a\phi_p,\psi_p,dx_p)}
 {|\eps(\tfrac12,\Sym^a\phi_p,\psi_p,dx_p)|},
\end{equation*}
where $a(\cdot)$ is the Artin conductor exponent and $\psi_p,dx_p$ are the
standard local components of the global additive character and self-dual
measure.  These conventions agree with \cite{Deligne,DMW}.  Then
\begin{equation}\label{eq:conductor-root-assembly}
 f_p(F(E))=\sum_{a,b}c_F(a,b)f_{p,a},
 \qquad
 W_p(F(E))=\prod_{a,b}W_{p,a}^{c_F(a,b)}.
\end{equation}
In particular,
\begin{equation}\label{eq:global-conductor}
 N_{F(E)}=\prod_{p\mid N_E}p^{f_p(F(E))}.
\end{equation}
\end{proposition}

\begin{proof}
Artin conductors are additive and epsilon factors are multiplicative under
direct sums.  The twist by $\omega_p^{-b}$ is unramified and positive real.
It changes the variable $s$ but not the conductor exponent or the phase of
the root number.
\end{proof}

\begin{remark}\label{rem:level-not-enough}
The integer $N_E$ alone does not determine the factors at its prime divisors.
One also needs the split or nonsplit multiplicative type, the quadratic
character in the potentially multiplicative case, or the inertia and
Frobenius action in the potentially good case.  At $p=2$ and $p=3$, the wild
filtration or equivalent minimal-model data are also needed.  For a fixed
elliptic curve these data are obtained from Tate's algorithm.
\end{remark}

\appendix

\section{Complete algebraic classification}\label{app:algebra}

This appendix makes the algebraic bookkeeping explicit before any local
Euler factor is taken.  Its basic ingredients are standard: the
highest-weight classification and Clebsch--Gordan formula may be found in
\cite{FultonHarris}, the plethysm multiplicities are given by the classical
Cayley--Sylvester formula, and the exact large-parameter support uses
Pak--Panova \cite{PakPanova}.  We organize these known identities into a
recursive classification adapted to the present problem.  We do not claim
the underlying invariant-theory formulae as new.

Every displayed decomposition satisfies the following useful dimension
check:
\begin{align*}
 \dim(W_1\oplus W_2)&=\dim W_1+\dim W_2,&
 \dim(W_1\otimes W_2)&=\dim W_1\dim W_2,\\
 \dim\Sym^mW&=\binom{\dim W+m-1}{m}.&&
\end{align*}

\subsection{The irreducible atoms}

\begin{proposition}\label{prop:all-atoms}
Fix a polynomial degree $D$.  The irreducible representations that can occur
in degree $D$ are exactly
\begin{equation}\label{eq:degree-D-atoms}
 R_{D-2b,b}=\Sym^{D-2b}V\otimes\det(V)^b,
 \qquad
 0\leq b\leq\left\lfloor\frac D2\right\rfloor.
\end{equation}
Every representation in this list occurs as an irreducible constituent of an
admissible expression.
\end{proposition}

\begin{proof}
The classification of polynomial representations of $\GL_2$ gives the first
assertion.  For the second, the Clebsch--Gordan formula gives
\begin{equation*}
 \Sym^{D-b}V\otimes\Sym^bV
 =\bigoplus_{j=0}^{b}R_{D-2j,j}.
\end{equation*}
The term with $j=b$ is the required representation.
\end{proof}

Thus the complete list across all degrees is
\begin{equation}\label{eq:all-Rab}
 \{R_{a,b}:a,b\geq0\}.
\end{equation}

\subsection{Tensor products}

\begin{proposition}[Clebsch--Gordan]\label{prop:CG}
For $a,c,u,v\geq0$,
\begin{equation}\label{eq:CG}
 R_{a,u}\otimes R_{c,v}
 =
 \bigoplus_{j=0}^{\min(a,c)}
 R_{a+c-2j,u+v+j}.
\end{equation}
\end{proposition}

\begin{proof}
The usual $\GL_2$ Clebsch--Gordan decomposition is
\begin{equation*}
 \Sym^aV\otimes\Sym^cV
 =\bigoplus_{j=0}^{\min(a,c)}
 \Sym^{a+c-2j}V\otimes\det(V)^j;
\end{equation*}
see, for example, \cite[Section~11.2]{FultonHarris}.  Tensoring by
$\det(V)^{u+v}$ gives \eqref{eq:CG}.
\end{proof}

In particular, the determinant powers in \eqref{eq:CG} cannot be omitted.
For a multiple tensor product there is a closed formula.

\begin{lemma}[Interval convolution]\label{lem:interval-convolution}
Let $d_1,\ldots,d_s\geq0$, put $S=\sum_i d_i$ and $H=\max_i d_i$, and write
\begin{equation*}
 \prod_{i=1}^{s}(1+q+\cdots+q^{d_i})
 =\sum_{k=0}^{S}a_kq^k,\qquad a_{-1}=0.
\end{equation*}
For $0\leq k\leq\lfloor S/2\rfloor$,
\begin{equation}\label{eq:interval-convolution}
 a_k-a_{k-1}>0\quad\Longleftrightarrow\quad k\leq S-H.
\end{equation}
\end{lemma}

\begin{proof}
Delete factors with $d_i=0$ and induct on the number of remaining factors.
If none remains, then $S=H=0$, $a_0-a_{-1}=1$, and the assertion is
immediate.  It is also immediate for one remaining factor.  Order the
factors so that
$d_s=H$, put $R=S-H$, and write
\begin{equation*}
 \prod_{i<s}(1+q+\cdots+q^{d_i})=\sum_{j=0}^{R}c_jq^j,
\end{equation*}
extending $c_j$ by zero outside $[0,R]$.  Convolution gives
\begin{equation}\label{eq:interval-difference}
 a_k-a_{k-1}=c_k-c_{k-H-1}.
\end{equation}
The sequence $(c_j)$ is palindromic and positive on $[0,R]$.  Let
$G=\max_{i<s}d_i\leq H$.  If
$0\leq k\leq\min(R,\lfloor(R+H)/2\rfloor)$ and $k\leq H$, the second term
in \eqref{eq:interval-difference} vanishes.  If $k>H$ and $k\leq R/2$, the
induction hypothesis makes $(c_j)$ strictly increasing between the two
indices in \eqref{eq:interval-difference}.  If $k>R/2$, palindromicity gives
$c_k=c_{R-k}$, while
\begin{equation*}
 k-H-1<R-k<R-G;
\end{equation*}
the induction hypothesis again gives $c_{k-H-1}<c_{R-k}$.  Thus the
difference is positive throughout the stated range.  If
$R<k\leq\lfloor(R+H)/2\rfloor$, then $H>R$ and both terms in
\eqref{eq:interval-difference} vanish.  This proves
\eqref{eq:interval-convolution}.
\end{proof}

\begin{proposition}[Multiple tensor products]\label{prop:multiple-tensor}
Let $r\geq1$ and $n_i,b_i\geq0$.  Put
\begin{equation*}
 M=\sum_{i=1}^{r}n_i,
 \qquad B=\sum_{i=1}^{r}b_i,
 \qquad
 A_k=[q^k]\prod_{i=1}^{r}(1+q+\cdots+q^{n_i}),
\end{equation*}
and put $A_{-1}=0$.  Then
\begin{equation}\label{eq:multiple-tensor}
 \bigotimes_{i=1}^{r}R_{n_i,b_i}
 =
 \bigoplus_{0\leq k\leq\lfloor M/2\rfloor}
 R_{M-2k,B+k}^{\oplus(A_k-A_{k-1})}.
\end{equation}
If $L=\max_i n_i$, then $A_k-A_{k-1}>0$ exactly for
\begin{equation}\label{eq:tensor-support}
 0\leq k\leq
 \min\left(M-L,\left\lfloor\frac M2\right\rfloor\right).
\end{equation}
\end{proposition}

\begin{proof}
On the diagonal torus $\operatorname{diag}(z,z^{-1})$, the weight $M-2k$
occurs with multiplicity $A_k$.  A summand $R_{M-2j,B+j}$ contributes once
to each weight $M-2k$ with $j\leq k\leq M-j$.  Hence, on the lower half,
its highest-weight multiplicity is $A_k-A_{k-1}$, proving
\eqref{eq:multiple-tensor}.  Lemma~\ref{lem:interval-convolution}, with
$S=M$ and $H=L$, proves \eqref{eq:tensor-support}.
\end{proof}

For example, a Tate line occurs in $\bigotimes_i\Sym^{n_i}V$ if and only if
$M$ is even and $\max_i n_i\leq M/2$.  Its multiplicity is
\begin{equation}\label{eq:tensor-Tate-mult}
 A_{M/2}-A_{M/2-1}.
\end{equation}
When all $n_i=0$, this records the trivial line $R_{0,0}$.

\subsection{Binary plethysm}

Let $p_{m,n}(r)$ be the number of partitions of $r$ whose Ferrers diagram fits
in an $m$ by $n$ rectangle.  Thus
\begin{equation}\label{eq:q-binomial}
 \sum_{r\geq0}p_{m,n}(r)q^r
 =\qbinom{m+n}{m}_q
 =\prod_{j=1}^{m}\frac{1-q^{n+j}}{1-q^j}.
\end{equation}
Put $p_{m,n}(-1)=0$ and
\begin{equation}\label{eq:mu-def}
 \mu_{m,n}(r)=p_{m,n}(r)-p_{m,n}(r-1).
\end{equation}
When $mn=0$, we extend the notation by
\begin{equation}\label{eq:mu-zero-convention}
 \mu_{m,n}(0)=1,\qquad \mu_{m,n}(r)=0\quad(r>0).
\end{equation}
This is the identity $\Sym^0(W)=\Sym^m(1)=1$.

\begin{theorem}[Cayley--Sylvester]\label{thm:Cayley-Sylvester}
For $m,n\geq1$ and $c\geq0$,
\begin{equation}\label{eq:Cayley-Sylvester}
 \Sym^m(R_{n,c})
 =
 \bigoplus_{0\leq r\leq\lfloor mn/2\rfloor}
 R_{mn-2r,mc+r}^{\oplus\mu_{m,n}(r)}.
\end{equation}
\end{theorem}

\begin{proof}
The weights of $\Sym^nV$ are $x^{n-j}y^j$, $0\leq j\leq n$.  A weight in
its $m$th symmetric power is a multiset of $m$ such weights.  After removing
the highest monomial $x^{mn}$, the exponent loss $r$ is counted by partitions
whose diagram fits in an $m$ by $n$ rectangle, hence has multiplicity
$p_{m,n}(r)$.  For $\GL_2$, subtracting the adjacent weight multiplicity
$p_{m,n}(r-1)$ extracts the highest-weight multiplicity.  The highest weight
$(mn-r,r)$ is $R_{mn-2r,r}$; the initial determinant twist contributes
$\det(V)^{mc}$.  This proves \eqref{eq:Cayley-Sylvester}.
\end{proof}

The multiplicities in \eqref{eq:Cayley-Sylvester} are not all one.  Some of
them vanish.  The next proposition lists every vanishing multiplicity.  Put
\begin{equation*}
 Z_{m,n}=\left\{0\leq r\leq\left\lfloor\frac{mn}{2}\right\rfloor:
 \mu_{m,n}(r)=0\right\}.
\end{equation*}

\begin{proposition}[Exact plethysm support]\label{prop:plethysm-support}
Assume $m,n\geq1$.  Let $\ell=\min(m,n)$ and $h=\max(m,n)$.  The set
$Z_{m,n}$ is given by the
following table.
\begin{center}
\begin{tabular}{@{}ll@{}}
\toprule
Condition & $Z_{m,n}$ \\
\midrule
$\ell=1$ & $\{1,2,\ldots,\lfloor h/2\rfloor\}$ \\
$\ell=2$ & odd $r$ with $1\leq r\leq h$ \\
$\ell=3$, $h$ odd & $\{1,\lfloor3h/2\rfloor\}$ \\
$\ell=3$, $h\equiv0\pmod4$ & $\{1,3h/2-1\}$ \\
$\ell=3$, $h\equiv2\pmod4$ & $\{1,3h/2-2,3h/2\}$ \\
$\ell=4$, $h=4$ & $\{1,5,7\}$ \\
$\ell=4$, $h\geq5$ & $\{1,2h-1\}$ \\
$\ell\geq5$, outside the exceptions below & $\{1\}$ \\
\bottomrule
\end{tabular}
\end{center}
For
\begin{equation}\label{eq:PP-exceptions}
 (\ell,h)\in
 \{(5,6),(5,10),(5,14),(6,7),(6,9),(6,11),(6,13),(7,10)\},
\end{equation}
one must add $\ell h/2$ to $Z_{m,n}$.  The remaining exception is
\begin{equation*}
 Z_{6,6}=\{1,17\}.
\end{equation*}
Every index outside this list has strictly positive multiplicity.
\end{proposition}

\begin{proof}
The multiplicity $\mu_{m,n}(r)$ is the coefficient of $q^r$ in
$(1-q)\qbinom{m+n}{m}_q$.  We spell out both parts of the support
calculation.

For $\ell=1$, the Gaussian polynomial is $1+q+\cdots+q^h$, so its
coefficient differences vanish at every $1\leq r\leq\lfloor h/2\rfloor$.

For $\ell=2$ and $0\leq r\leq h$,
\begin{equation*}
 (1-q)\qbinom{h+2}{2}_q
 =\frac{(1-q^{h+1})(1-q^{h+2})}{1-q^2},
\end{equation*}
so the coefficient is $1$ for even $r$ and $0$ for odd $r$.  For
$\ell=3$, put
\begin{equation*}
 A(t)=\#\{(x,y)\in\Z_{\geq0}^2:2x+3y=t\},\qquad A(t)=0\ (t<0).
\end{equation*}
For $r\leq\lfloor3h/2\rfloor$, two shifted numerator terms cannot occur, and
the coefficient is
\begin{equation}\label{eq:l3-coeff-check}
 A(r)-\sum_{i=1}^{3}A(r-h-i),
 \qquad
 A(t)=
 \begin{cases}
  \lfloor t/6\rfloor,&t\equiv1\pmod6,\\
  \lfloor t/6\rfloor+1,&t\not\equiv1\pmod6.
 \end{cases}\qquad(t\geq0).
\end{equation}
Substitution in the six residue classes gives, besides $r=1$, respectively
the endpoint zeros
\begin{equation*}
 \lfloor3h/2\rfloor\ (h\text{ odd}),\qquad
 3h/2-1\ (h\equiv0\!\!\pmod4),\qquad
 3h/2-2,\ 3h/2\ (h\equiv2\!\!\pmod4),
\end{equation*}
and positive coefficients elsewhere.

For $\ell=4$, put
\begin{equation*}
 B(t)=\#\{(x,y,z)\in\Z_{\geq0}^3:2x+3y+4z=t\}
      =\sum_{z=0}^{\lfloor t/4\rfloor}A(t-4z),
\end{equation*}
again taking $B(t)=0$ for $t<0$.  Since two numerator shifts exceed $2h$,
the required coefficient for $r\leq2h$ is
\begin{equation}\label{eq:l4-coeff-check}
 B(r)-\sum_{i=1}^{4}B(r-h-i).
\end{equation}
For an entirely explicit check, let
\begin{equation*}
 \rho(v)=[q^v]\frac1{(1-q)(1-q^2)(1-q^3)},\qquad \rho(v)=0\ (v<0).
\end{equation*}
The elementary formula for partitions into parts at most three gives
\begin{equation}\label{eq:rho-explicit}
 \rho(v)=\left\lfloor\frac{v^2+6v+12}{12}\right\rfloor\quad(v\geq0),
 \qquad
 B(2u)=\rho(u),\quad B(2u+1)=\rho(u-1).
\end{equation}
Substitution in \eqref{eq:l4-coeff-check}, separated by parity and then
modulo $6$, gives zeros $1,5,7$ for $h=4$, and $1,2h-1$ for $h\geq5$;
all other coefficients are positive.  This proves every row with
$\ell\leq4$ rather than merely appealing to unimodality.

For $\ell\geq5$, Theorem~6 and the corrected exceptional list of Pak and
Panova give strict increase through the midpoint except at $r=1$ and the
nine plateaux below \cite[Theorem~6 and the correction]{PakPanova}:
the displayed integer checks follow recursively from
\begin{equation*}
 \qbinom{\ell+h}{\ell}_q
 =\qbinom{\ell+h-1}{\ell}_q
 +q^h\qbinom{\ell+h-1}{\ell-1}_q.
\end{equation*}
\begin{center}
\begin{tabular}{@{}ccc@{}}
\toprule
$(\ell,h)$ & missing $r$ & equal adjacent coefficients \\
\midrule
$(5,6)$ & $15$ & $p(14)=p(15)=32$ \\
$(5,10)$ & $25$ & $p(24)=p(25)=141$ \\
$(5,14)$ & $35$ & $p(34)=p(35)=414$ \\
$(6,6)$ & $17$ & $p(16)=p(17)=55$ \\
$(6,7)$ & $21$ & $p(20)=p(21)=94$ \\
$(6,9)$ & $27$ & $p(26)=p(27)=227$ \\
$(6,11)$ & $33$ & $p(32)=p(33)=480$ \\
$(6,13)$ & $39$ & $p(38)=p(39)=920$ \\
$(7,10)$ & $35$ & $p(34)=p(35)=734$ \\
\bottomrule
\end{tabular}
\end{center}
Here $p(r)=p_{\ell,h}(r)$.  In eight cases the missing index is the
midpoint $\ell h/2$; the exceptional $(6,6)$ plateau is at $17$, not at the
midpoint $18$.  Together these calculations give exactly the displayed
support.  Finally, $r=1$ always vanishes when $m,n\geq2$, since
$p_{m,n}(0)=p_{m,n}(1)=1$.
\end{proof}

Combining Theorem~\ref{thm:Cayley-Sylvester} with
Proposition~\ref{prop:plethysm-support} gives the nonredundant formula
\begin{equation}\label{eq:plethysm-nonredundant}
 \Sym^m(R_{n,c})
 =
 \bigoplus_{\substack{0\leq r\leq\lfloor mn/2\rfloor\\r\notin Z_{m,n}}}
 R_{mn-2r,mc+r}^{\oplus\mu_{m,n}(r)}.
\end{equation}

\subsection{Direct sums and arbitrary nesting}

For representations $W_1,\ldots,W_t$,
\begin{equation}\label{eq:sym-direct-sum}
 \Sym^m\left(\bigoplus_{i=1}^{t}W_i\right)
 =
 \bigoplus_{\substack{m_1+\cdots+m_t=m\\m_i\geq0}}
 \bigotimes_{i=1}^{t}\Sym^{m_i}W_i.
\end{equation}
Equations \eqref{eq:multiple-tensor}, \eqref{eq:Cayley-Sylvester}, and
\eqref{eq:sym-direct-sum} calculate every finite nesting.  The following
character formula is an equivalent one-line algorithm.  If
\begin{equation*}
 W=\bigoplus_{a,b}R_{a,b}^{\oplus c(a,b)},
\end{equation*}
then
\begin{equation}\label{eq:character-generating}
 \sum_{m\geq0}\ch(\Sym^mW)t^m
 =
 \prod_{a,b}\prod_{j=0}^{a}
 \left(1-tx^{a+b-j}y^{b+j}\right)^{-c(a,b)}.
\end{equation}
If $H(x,y)$ is a symmetric character, the multiplicity of $R_{a,b}$ in it is
\begin{equation}\label{eq:character-extraction}
 [x^{a+b}y^b]H-[x^{a+b+1}y^{b-1}]H,
\end{equation}
where the second coefficient is zero for $b=0$.

For clarity, we now distinguish an irreducible constituent from an entire
expression.  For example, $\det V$ occurs in
\begin{equation*}
 V\otimes V=\Sym^2V\oplus\det V,
\end{equation*}
but it cannot be obtained as an entire expression unless taking a direct
summand is added as an allowed operation.

\begin{proposition}[Exact set of whole expressions]\label{prop:whole-expressions}
Fix $D\geq1$.  Write a homogeneous representation of degree $D$ as
\begin{equation*}
 \bigoplus_{b=0}^{\lfloor D/2\rfloor}
 R_{D-2b,b}^{\oplus c_b}
\end{equation*}
and identify it with its multiplicity vector
$c=(c_0,\ldots,c_{\lfloor D/2\rfloor})$.

Let $\mathcal C_D$ be the finite set of vectors obtained from connected
operation trees.  Their leaves are $\Sym^nV$, and their internal vertices are
tensor products or symmetric powers.  Let $\star$ denote the product induced
by \eqref{eq:CG}, and let $\sigma_m$ denote the operation induced by
\eqref{eq:character-generating}.  Then
\begin{equation}\label{eq:CD-recursion}
 \mathcal C_D
 =
 \{(1,0,\ldots,0)\}
 \cup
 \bigcup_{i=1}^{D-1}
 \left(\mathcal C_i\star\mathcal C_{D-i}\right)
 \cup
 \bigcup_{\substack{m\mid D\\m\geq2}}
 \sigma_m\left(\mathcal C_{D/m}\right).
\end{equation}
After adjoining the zero object for bookkeeping, the set of all degree $D$
outputs is the finitely generated submonoid
\begin{equation}\label{eq:MD}
 \mathcal M_D=\Z_{\geq0}\mathcal C_D.
\end{equation}
Every general, possibly mixed-degree output is a finite direct sum of elements
of the monoids $\mathcal M_D$.
\end{proposition}

\begin{proof}
Use distributivity and \eqref{eq:sym-direct-sum} to move every direct sum to
the outermost level.  Each remaining summand is a connected operation tree.
Its root is either a leaf, a tensor product, or a symmetric power.  These are
exactly the three terms in \eqref{eq:CD-recursion}.  Conversely, every vector
in the recursion is produced by an allowed operation.  Taking outer direct
sums gives precisely the nonnegative span in \eqref{eq:MD}.
\end{proof}

In degree zero we put
\begin{equation}\label{eq:degree-zero-monoid}
 \mathcal C_0=\{(1)\},\qquad
 \mathcal M_0=\Z_{\geq0}(1).
\end{equation}
This records finite sums of the trivial representation.  Vertices labelled
$\Sym^1$ and tensor products with $1$ may be removed from a connected tree;
$\Sym^0(W)=1$ is already accounted for by \eqref{eq:degree-zero-monoid}.
Thus the recursion for $D\geq1$ remains finite and complete.
To compute it, construct $\mathcal C_1,\ldots,\mathcal C_D$ in increasing
degree from \eqref{eq:CD-recursion}, deduplicating the finitely many vectors
after each union, and then take the nonnegative span for $\mathcal M_D$.

Every vector in $\mathcal C_D$ has first coordinate $1$.  Thus a vector
$c=(c_0,\ldots,c_{\lfloor D/2\rfloor})$ is realizable if and only if it is a
sum of exactly $c_0$ vectors from $\mathcal C_D$.  The first five connected
sets are as follows.
\begin{align*}
 \mathcal C_1={}&\{(1)\},\\
 \mathcal C_2={}&\{(1,0),(1,1)\},\\
 \mathcal C_3={}&\{(1,0),(1,1),(1,2)\},\\
 \mathcal C_4={}&\{(1,0,0),(1,0,1),(1,1,0),(1,1,1),\\
 &\hspace{24mm}(1,1,2),(1,2,1),(1,3,2)\},\\
 \mathcal C_5={}&\{(1,0,0),(1,1,0),(1,1,1),(1,2,1),\\
 &\hspace{24mm}(1,2,2),(1,2,3),(1,3,3),(1,4,5)\}.
\end{align*}
The coordinates correspond in order to
\begin{equation*}
 R_{D,0},R_{D-2,1},R_{D-4,2},\ldots.
\end{equation*}

\subsection{The two families written in the original argument}

The first family is
\begin{equation*}
 \Sym^m\left(\bigoplus_{i=1}^{t}\Sym^{n_i}V\right).
\end{equation*}
Choose $m_1+\cdots+m_t=m$.  For every $i$, choose
\begin{equation*}
 0\leq r_i\leq\left\lfloor\frac{m_in_i}{2}\right\rfloor,
 \qquad
 \ell_i=m_in_i-2r_i,
\end{equation*}
with multiplicity $\mu_{m_i,n_i}(r_i)$.  Put
\begin{equation*}
 A_{\boldsymbol\ell}(k)
 =[q^k]\prod_{i=1}^{t}(1+q+\cdots+q^{\ell_i}).
\end{equation*}
Here, if $m_in_i=0$, the convention is $r_i=\ell_i=0$ and the multiplicity
is $1$; no other choice is allowed for that index.
Then the full decomposition is
\begin{equation}\label{eq:first-original-family}
 \bigoplus_{\substack{m_1+\cdots+m_t=m\\
 0\leq r_i\leq\lfloor m_in_i/2\rfloor\\
 0\leq k\leq\lfloor\sum_i\ell_i/2\rfloor}}
 R_{\sum_i\ell_i-2k,\ \sum_i r_i+k}^{\oplus M(\boldsymbol m,\boldsymbol r,k)},
\end{equation}
where
\begin{equation}\label{eq:first-family-mult}
 M(\boldsymbol m,\boldsymbol r,k)
 =
 \left(\prod_i\mu_{m_i,n_i}(r_i)\right)
 \left(A_{\boldsymbol\ell}(k)-A_{\boldsymbol\ell}(k-1)\right).
\end{equation}
Zero multiplicities are omitted, and equal terms are combined.

For the second family one has the simpler complete formula
\begin{equation}\label{eq:second-original-family}
 \left(\bigoplus_{i=1}^{u}\Sym^{n_i}V\right)
 \otimes
 \left(\bigoplus_{j=1}^{v}\Sym^{m_j}V\right)
 =
 \bigoplus_{i=1}^{u}\bigoplus_{j=1}^{v}
 \bigoplus_{k=0}^{\min(n_i,m_j)}
 R_{n_i+m_j-2k,k}.
\end{equation}
The upper bound is $\min(n_i,m_j)$, not $\min(i,j)$.

\subsection{Low-degree decompositions}

\begin{center}
\renewcommand{\arraystretch}{1.15}
\begin{tabular}{@{}ll@{}}
\toprule
Expression & Decomposition \\
\midrule
$R_{1,0}\otimes R_{1,0}$ & $R_{2,0}\oplus R_{0,1}$ \\
$R_{1,0}\otimes R_{2,0}$ & $R_{3,0}\oplus R_{1,1}$ \\
$R_{2,0}\otimes R_{2,0}$ & $R_{4,0}\oplus R_{2,1}\oplus R_{0,2}$ \\
$R_{2,0}\otimes R_{3,0}$ & $R_{5,0}\oplus R_{3,1}\oplus R_{1,2}$ \\
$R_{3,0}\otimes R_{3,0}$ &
$R_{6,0}\oplus R_{4,1}\oplus R_{2,2}\oplus R_{0,3}$ \\
$\Sym^2(R_{2,0})$ & $R_{4,0}\oplus R_{0,2}$ \\
$\Sym^3(R_{2,0})$ & $R_{6,0}\oplus R_{2,2}$ \\
$\Sym^4(R_{2,0})$ & $R_{8,0}\oplus R_{4,2}\oplus R_{0,4}$ \\
$\Sym^2(R_{3,0})$ & $R_{6,0}\oplus R_{2,2}$ \\
$\Sym^3(R_{3,0})$ & $R_{9,0}\oplus R_{5,2}\oplus R_{3,3}$ \\
$\Sym^4(R_{3,0})$ &
$R_{12,0}\oplus R_{8,2}\oplus R_{6,3}\oplus R_{4,4}\oplus R_{0,6}$ \\
\bottomrule
\end{tabular}
\end{center}

As a complete nested example,
\begin{equation}\label{eq:nested-running-example}
 \Sym^2(V\otimes V)
 =R_{4,0}\oplus R_{2,1}\oplus R_{0,2}^{\oplus2}.
\end{equation}
Consequently, for every elliptic curve and every prime, good or bad,
\begin{equation}\label{eq:nested-running-L}
 L_p(s,\Sym^2(E\otimes E))
 =L_p(s,\Sym^4E)L_p(s-1,\Sym^2E)\zeta_p(s-2)^2.
\end{equation}
At a bad prime each factor on the right is read from the atomic local table;
one does not apply $\Sym^2$ to the visible roots of $L_p(s,E\otimes E)$.

\subsection{Tate and zeta factors}

For a non-CM elliptic curve, a motivic zeta factor occurs exactly when
$c_F(0,b)>0$.  The corresponding factor is
\begin{equation}\label{eq:zeta-factor}
 L(s-b,\Sym^0E)=\zeta(s-b).
\end{equation}
For two factors, $\Sym^aV\otimes\Sym^cV$ contains a Tate line if and only if
$a=c$.  It is $\det(V)^a$ and contributes $\zeta(s-a)$.  For
$\Sym^m(\Sym^nV)$, a Tate line occurs precisely when $mn$ is even and
\begin{equation*}
 \mu_{m,n}(mn/2)>0.
\end{equation*}
Its multiplicity is this partition difference and its factor is
$\zeta(s-mn/2)$.

\section{Finite local factors for elliptic curves over \texorpdfstring{$\Q$}{Q}}
\label{app:local}

We now list the local input needed in \eqref{eq:P-F-assembly}.  These atomic
Euler factors and root numbers are reorganized from the calculations of
Dummigan--Martin--Watkins \cite{DMW}.  The wild conductor and minimal-model
data at $2$ and $3$ are taken from Martin--Watkins
\cite{MartinWatkins}.  Our contribution here is the translation into the
$R_{a,b}$ notation and the assembly for arbitrary admissible expressions.
No new atomic local calculation is claimed.  The base field is always
$\Q_p$.  Some formulae change over a general local field with even residue
degree.

Here is the route through the appendix for a fixed minimal Weierstrass model.
First determine the reduction type with Tate's algorithm.  If
$v_p(j_E)<0$, use the multiplicative or potentially multiplicative rows.  If
$v_p(j_E)\geq0$ and $p\geq5$, use the Kodaira table below.  At $p=3$, use
$f_{3,1}$ and $v_3(\Delta)$; at $p=2$, first choose the minimal quadratic
twist specified below.  The resulting $P_{p,n}$, $f_{p,n}$, and $W_{p,n}$
are finally inserted into \eqref{eq:P-F-assembly} and
\eqref{eq:conductor-root-assembly}.

For potentially good reduction, the monodromy operator is zero and the
normalized representation
\begin{equation}\label{eq:finite-WD-image}
 \rho_p^u=r_p\otimes\omega_p^{1/2}
\end{equation}
has finite inertia image.  In the tables, $I=\rho_p^u(I_p)$ denotes that
finite image.  We use $G$ only as shorthand for the group generated by $I$
and normalized Frobenius; $G$ need not be finite.  Saying that $G$ is
nonabelian means that Frobenius acts nontrivially on $I$, not that a cyclic
inertia group is itself nonabelian.

In a potentially good case, let $\gamma_p,\delta_p$ be the two eigenvalues of
$r_p(\Phi_p)$ on the complete two-dimensional Weil--Deligne space; then
$\gamma_p\delta_p=p$.  They need not occur among the visible roots of
$P_{p,E}(T)$.  They are determined from the full local Weil--Deligne data (equivalently,
by the good-reduction model over a suitable extension), as in
\cite[Proposition~5.1]{DMW}.

For navigation: when $I$ is cyclic, ``$G$ abelian'' means that Frobenius
commutes with $I$ and one uses item~5 below; ``$G$ nonabelian'' means that
Frobenius acts by inversion and one uses item~6.  For noncyclic inertia, use
items~7--9 directly.

\subsection{Inertia invariants}

For an inertia group $I$, put
\begin{equation}\label{eq:beta-epsilon}
 \beta_n(I)=\dim(\Sym^nV)^I,
 \qquad
 \epsilon_n(I)=n+1-\beta_n(I).
\end{equation}
When $N=0$, $\beta_n(I)$ is the degree of $P_{p,n}$, while
$\epsilon_n(I)$ is the tame contribution to the conductor exponent.
For $e=2,3,4,6$,
\begin{equation}\label{eq:beta-cyclic}
 \beta_n(C_e)=
 \#\{0\leq j\leq n:e\mid n-2j\}.
\end{equation}
For each of the three noncyclic inertia groups, $\beta_n(I)=0$ when $n$ is
odd.  For even $n$, define
\begin{equation*}
 u_n=
 \begin{cases}
 1,&n\equiv0\pmod6,\\
 0,&n\equiv2\pmod6,\\
 -1,&n\equiv4\pmod6.
 \end{cases}
\end{equation*}
Then
\begin{align}
 \beta_n(Q_8)
 &=\frac{n+1+3(-1)^{n/2}}4,\label{eq:beta-Q8}\\
 \beta_n(C_3\rtimes C_4)
 &=\frac{n+1+2u_n+3(-1)^{n/2}}6,\label{eq:beta-G24}\\
 \beta_n(\SL_2(\mathbb F_3))
 &=\frac{n+1+8u_n+3(-1)^{n/2}}{12}.
 \label{eq:beta-SL2}
\end{align}
These character calculations agree with
\cite[Table~1]{MartinWatkins}.  The complementary number
$\epsilon_n(I)$ is the tame part of the Artin conductor; despite the letter,
it is not an epsilon factor.

\subsection{Atomic Euler factors by reduction type}

The following nine cases exhaust the reduction, inertia, and Frobenius-action
families needed to compute the atomic Euler polynomials of an elliptic curve
over $\Q_p$, after retaining the indicated Frobenius and minimal-model data.
Indeed, nonzero monodromy gives the multiplicative or
potentially multiplicative cases.  When $N=0$, the finite inertia image is
trivial, cyclic of order $2,3,4,$ or $6$, or one of
$C_3\rtimes C_4,Q_8,\SL_2(\mathbb F_3)$; see
\cite[Section~4]{DMW}.  We first give the atomic polynomial $P_{p,n}(T)$ defined in
\eqref{eq:Ppa}.  For $R_{n,b}$, replace $T$ by $p^bT$.

\begin{enumerate}
\item If $E$ has good reduction, then
\begin{equation}\label{eq:local-good}
 P_{p,n}(T)=
 \prod_{j=0}^{n}
 (1-\alpha_p^{n-j}\beta_p^jT).
\end{equation}

\item If $E$ has split multiplicative reduction, then
\begin{equation}\label{eq:local-split}
 P_{p,n}(T)=1-T.
\end{equation}

\item If $E$ has nonsplit multiplicative reduction, then
\begin{equation}\label{eq:local-nonsplit}
 P_{p,n}(T)=1-(-1)^nT.
\end{equation}

\item If $E$ has additive potentially multiplicative reduction, then
\begin{equation}\label{eq:local-add-pot-mult}
 P_{p,n}(T)=
 \begin{cases}
 1,&n\ \text{odd},\\
 1-T,&n\ \text{even}.
 \end{cases}
\end{equation}
These first four cases are \cite[Proposition~4.1]{DMW}, together with the
good-reduction definition.

\item Suppose $E$ has potentially good reduction, $I=C_e$ with
$e\in\{2,3,4,6\}$, and Frobenius commutes with inertia.  For suitable
$\gamma_p,\delta_p$ as defined above,
\begin{equation}\label{eq:local-cyclic-commuting}
 P_{p,n}(T)=
 \prod_{\substack{0\leq j\leq n\\e\mid n-2j}}
 (1-\gamma_p^{n-j}\delta_p^jT).
\end{equation}
This is \cite[Proposition~5.1]{DMW}.

\item Suppose $I=C_e$, with $e\in\{3,4,6\}$, and Frobenius inverts inertia.
If $n$ is odd, then
\begin{equation}\label{eq:local-cyclic-inverting-odd}
 P_{p,n}(T)=
 (1-(-p)^nT^2)^{\beta_n(C_e)/2}.
\end{equation}
If $n$ is even and $x=(-p)^{n/2}T$, then
\begin{equation}\label{eq:local-cyclic-inverting-even}
 P_{p,n}(T)=
 (1-x)^{\lceil\beta_n(C_e)/2\rceil}
 (1+x)^{\lfloor\beta_n(C_e)/2\rfloor}.
\end{equation}
These two formulas are \cite[Proposition~6.1]{DMW}.

\item If $p=3$ and $I=C_3\rtimes C_4$, then $P_{p,n}(T)=1$ for odd $n$.
For even $n$, formula \eqref{eq:local-cyclic-inverting-even} holds with
$\beta_n(C_3\rtimes C_4)$ \cite[Section~11]{DMW}.

\item If $p=2$ and $I=Q_8$, then $P_{p,n}(T)=1$ for odd $n$.  For even $n$,
formula \eqref{eq:local-cyclic-inverting-even} holds with $\beta_n(Q_8)$
\cite[Proposition~8.1]{DMW}.

\item If $p=2$ and $I=\SL_2(\mathbb F_3)$, then $P_{p,n}(T)=1$ for odd $n$.
For even $n$, formula \eqref{eq:local-cyclic-inverting-even} holds with
$\beta_n(\SL_2(\mathbb F_3))$ \cite[Section~10]{DMW}.
\end{enumerate}

The last two formulae use the fact that the residue degree of $\Q_p$ is one.
Over a local field of even residue degree the $\SL_2(\mathbb F_3)$ case can
also contain cubic roots of unity.

\subsection{Direct ramified formulae for the composite operations}

Let
\begin{equation*}
 T_0=\bigotimes_iR_{n_i,b_i},\qquad
 M=\sum_i n_i,\qquad B=\sum_i b_i,
\end{equation*}
and let $A_k-A_{k-1}$ be the multiplicities from
Proposition~\ref{prop:multiple-tensor}.  At every prime,
\begin{equation}\label{eq:ramified-tensor-master}
 L_p(s,T_0(E))
 =
 \prod_k
 L_p(s-B-k,\Sym^{M-2k}E)^{A_k-A_{k-1}}.
\end{equation}
For split multiplicative reduction this becomes
\begin{equation}\label{eq:ramified-tensor-split}
 \prod_k(1-p^{B+k-s})^{-(A_k-A_{k-1})}.
\end{equation}
For nonsplit multiplicative reduction it becomes
\begin{equation}\label{eq:ramified-tensor-nonsplit}
 \prod_k
 (1-(-1)^Mp^{B+k-s})^{-(A_k-A_{k-1})}.
\end{equation}
For additive potentially multiplicative reduction it equals $1$ when $M$ is
odd, and it equals \eqref{eq:ramified-tensor-split} when $M$ is even.

For a binary plethysm, put $Q=mn$.  Then
\begin{equation}\label{eq:ramified-plethysm-master}
 L_p(s,\Sym^m(R_{n,c})(E))
 =
 \prod_{r=0}^{\lfloor Q/2\rfloor}
 L_p(s-mc-r,\Sym^{Q-2r}E)^{\mu_{m,n}(r)}.
\end{equation}
Again, every atomic factor on the right side is taken from the complete local
type and not from the roots of $L_p(s,\Sym^nE)$.

\subsection{Multiplicative and potentially multiplicative conductors}

In the potentially multiplicative case there is a quadratic character
$\eta_p$ such that
\begin{equation}\label{eq:pot-mult-WD}
 \tau_p\simeq\eta_p\omega_p^{-1}\otimes\spc(2),
 \qquad
 \Sym^n\tau_p\simeq
 \eta_p^n\omega_p^{-n}\otimes\spc(n+1).
\end{equation}
Here $\spc(m)$ is the special Weil--Deligne representation with basis
$e_0,\ldots,e_{m-1}$, a single Jordan block
$Ne_j=e_{j+1}$, $Ne_{m-1}=0$, and Weil action
$r(g)e_j=\omega_p(g)^je_j$ (up to the common unramified twist displayed in
\eqref{eq:pot-mult-WD}).
Write $w_p=W_p(E)$.  The complete conductor and root-number table for these
cases is
\begin{center}
\renewcommand{\arraystretch}{1.2}
\begin{tabular}{@{}p{.26\textwidth}p{.31\textwidth}p{.28\textwidth}@{}}
\toprule
Reduction type & $f_{p,n}$ & $W_{p,n}$ \\
\midrule
Good & $0$ & $1$ \\
Split multiplicative & $n$ & $(-1)^n$ \\
Nonsplit multiplicative & $n$ & $1$ \\
Additive potentially multiplicative, $n$ odd
& $(n+1)a(\eta_p)$ & $w_p^{(n+1)/2}$ \\
Additive potentially multiplicative, $n$ even
& $n$ & $1$ \\
\bottomrule
\end{tabular}
\end{center}
For an odd prime, a ramified quadratic $\eta_p$ has $a(\eta_p)=1$ and
\begin{equation*}
 w_p=\left(\frac{-1}{p}\right).
\end{equation*}
At $p=2$, one has $a(\eta_2)=2$ or $3$.  The value $w_2$ requires the
specific quadratic character or minimal model.
This table is \cite[Proposition~4.1]{DMW}.

\begin{proposition}[Even symmetric powers]\label{prop:even-root}
For every finite prime $p$ and every even $n$,
\begin{equation}\label{eq:even-root}
 W_{p,n}=1.
\end{equation}
\end{proposition}

This is \cite[Proposition~7.1]{DMW}.  The remaining tables therefore list
root numbers only for odd
$n$.

\subsection{Potentially good reduction for \texorpdfstring{$p\geq5$}{p at least 5}}

Let $\Delta$ be a minimal discriminant.  The local types are
\begin{center}
\renewcommand{\arraystretch}{1.15}
\begin{tabular}{@{}llll@{}}
\toprule
Kodaira type & $v_p(\Delta)$ & Reduction & $e=|I|$ \\
\midrule
$I_0$ & $0$ & good & $1$ \\
$I_m$ & $m\geq1$ & multiplicative & -- \\
$II$ & $2$ & potentially good & $6$ \\
$III$ & $3$ & potentially good & $4$ \\
$IV$ & $4$ & potentially good & $3$ \\
$I_0^*$ & $6$ & potentially good & $2$ \\
$IV^*$ & $8$ & potentially good & $3$ \\
$III^*$ & $9$ & potentially good & $4$ \\
$II^*$ & $10$ & potentially good & $6$ \\
$I_m^*$, $m\geq1$ & $m+6$ & potentially multiplicative & -- \\
\bottomrule
\end{tabular}
\end{center}
In the potentially good rows,
\begin{equation}\label{eq:e-Kodaira}
 e=\frac{12}{\gcd(12,v_p(\Delta))},
 \qquad
 f_{p,n}=\epsilon_n(C_e).
\end{equation}
Frobenius commutes with inertia if and only if
\begin{equation}\label{eq:Frob-commute}
 p\equiv1\pmod e.
\end{equation}
Otherwise it acts by inversion.  The basic root number is
\begin{equation}\label{eq:wp-tame}
 w_p=
 \begin{cases}
 \left(\dfrac{-1}{p}\right),&e=2\ \text{or}\ 6,\\[1mm]
 \left(\dfrac{-3}{p}\right),&e=3,\\[1mm]
 \left(\dfrac{-2}{p}\right),&e=4.
 \end{cases}
\end{equation}
For $e=2,4,6$ and odd $n$,
\begin{equation}\label{eq:root-e246}
 W_{p,n}=
 \begin{cases}
 w_p,&n\equiv1\pmod4,\\
 1,&n\equiv3\pmod4.
 \end{cases}
\end{equation}
For $e=3$ and $p\equiv1\pmod3$, all the root numbers are $1$.  For $e=3$
and $p\equiv2\pmod3$, the values for
$n\equiv1,3,5,7,9,11\pmod{12}$ are respectively
\begin{equation}\label{eq:e3-root-sequence}
 -1,\ -1,\ 1,\ -1,\ -1,\ 1.
\end{equation}
These tame formulas are \cite[Propositions~5.1 and 6.1]{DMW}.

\subsection{Potentially good reduction at \texorpdfstring{$p=3$}{p=3}}

Assume $v_3(j_E)\geq0$.  Put
\begin{equation*}
 \mathfrak f_3=f_{3,1}=v_3(N_E),\qquad d=v_3(\Delta).
\end{equation*}
If $\mathfrak f_3=0$, the curve has good reduction and item~1 of the atomic
list applies.  In the remainder of this subsection assume
$\mathfrak f_3>0$.
The possible inertia groups are
\begin{center}
\begin{tabular}{@{}lll@{}}
\toprule
Condition & $I$ & $G$ \\
\midrule
$\mathfrak f_3=2$, $d$ even & $C_2$ & abelian \\
$\mathfrak f_3=2$, $d$ odd & $C_4$ & nonabelian \\
$\mathfrak f_3=4$, $4\mid d$ & $C_3$ & either \\
$\mathfrak f_3=4$, $4\nmid d$ & $C_6$ & either \\
$\mathfrak f_3=3$ or $5$ & $C_3\rtimes C_4$ & nonabelian \\
\bottomrule
\end{tabular}
\end{center}
For $C_3$ and $C_6$, choose, as in Martin--Watkins, a local quadratic twist
$E^\flat/\Q_3$ having minimal conductor exponent, and take a minimal integral
Weierstrass model of it.  Put
\begin{equation*}
 \widehat c_4=c_4(E^\flat),\qquad
 \widehat c_6=c_6(E^\flat).
\end{equation*}
All congruence tests below refer to this Martin--Watkins choice.  If
$\widehat c_4\equiv9\pmod{27}$, the group is abelian if and only if
\begin{equation}\label{eq:p3-abelian-test1}
 \widehat c_6\equiv\pm108\pmod{243}.
\end{equation}
If $27\mid\widehat c_4$, it is abelian if and only if
\begin{equation}\label{eq:p3-abelian-test2}
 \widehat c_4\equiv27\pmod{81}.
\end{equation}

Write
\begin{equation*}
 f_{3,n}=\epsilon_n(I)+\delta_{3,n}.
\end{equation*}
The wild term is
\begin{center}
\begin{tabular}{@{}lll@{}}
\toprule
$I$ & $\mathfrak f_3$ & $\delta_{3,n}$ \\
\midrule
$C_2$ & $2$ & $0$ \\
$C_4$ & $2$ & $0$ \\
$C_3$ & $4$ & $\epsilon_n(C_3)$ \\
$C_6$ & $4$ & $\epsilon_n(C_3)$ \\
$C_3\rtimes C_4$ & $3$ & $\tfrac12\epsilon_n(C_3)$ \\
$C_3\rtimes C_4$ & $5$ & $\tfrac32\epsilon_n(C_3)$ \\
\bottomrule
\end{tabular}
\end{center}
These are the specializations to $\Q_3$ of
\cite[Tables~2 and 3]{MartinWatkins}.  Martin--Watkins tabulate the wild
exponent $\delta_{3,1}$; the column $\mathfrak f_3$ here is the total
conductor $\epsilon_1(I)+\delta_{3,1}$.

For odd $n$, the root numbers are given below.  The six entries correspond to
$n\equiv1,3,5,7,9,11\pmod{12}$.
\begin{center}
\begin{tabular}{@{}ll@{}}
\toprule
Local type & Root-number sequence \\
\midrule
$C_2$ & $-1,1,-1,1,-1,1$ \\
$C_4$ & $1,1,1,1,1,1$ \\
$C_3$ & $1,1,1,1,1,1$ \\
$C_6$, abelian & $-1,1,-1,1,-1,1$ \\
$C_6$, nonabelian & $1,-1,-1,-1,1,1$ \\
$C_3\rtimes C_4$ & $w_3,w_3,-1,w_3,w_3,1$ \\
\bottomrule
\end{tabular}
\end{center}
In the last row, $w_3=W_3(E)$ is determined by the minimal model and not by
$I$ and $\mathfrak f_3$ alone.
The cyclic rows follow from \cite[Propositions~5.1 and 6.1]{DMW}; the final
row is the $\Q_3$ case of \cite[Section~11]{DMW}.

\subsection{Potentially good reduction at \texorpdfstring{$p=2$}{p=2}}

Assume $v_2(j_E)\geq0$.  Among the local quadratic twists of $E/\Q_2$,
choose $E^\flat$ with minimal conductor exponent, using the choice convention
of \cite[Table~4]{MartinWatkins}, and take a minimal integral Weierstrass
model.  Put
\begin{equation*}
 \mathfrak f_2=f_{2,1}=v_2(N_E),\qquad
 m_\flat=v_2(N_{E^\flat}).
\end{equation*}
Good reduction is already covered by item~1; the tables below include its
appearance as the $m_\flat=0$, $v_2(N_E)=0$ subcase.
The possibilities are
\begin{center}
\begin{tabular}{@{}ll@{}}
\toprule
$m_\flat$ & Inertia group \\
\midrule
$0$ & $C_1$ if $v_2(N_E)=0$, otherwise $C_2$ \\
$2$ & $C_3$ if $v_2(N_E)=2$, otherwise $C_6$ \\
$3$ or $7$ & $\SL_2(\mathbb F_3)$ \\
$5$ & $Q_8$ \\
$8$ & $Q_8$ if $2^9\mid c_6(E^\flat)$, otherwise $C_4$ \\
\bottomrule
\end{tabular}
\end{center}
For $I=C_1,C_2$, the group $G$ is abelian.  For
$I=C_3,C_6,Q_8,\SL_2(\mathbb F_3)$, the group $G$ generated after Frobenius
is included is nonabelian; the statement does not make the cyclic groups
$C_3,C_6$ nonabelian.  For $I=C_4$, the group $G$ is abelian if
\begin{equation*}
 c_4(E^\flat)\equiv96\pmod{128},
\end{equation*}
and nonabelian if
\begin{equation*}
 c_4(E^\flat)\equiv32\pmod{128}.
\end{equation*}

Write
\begin{equation*}
 f_{2,n}=\epsilon_n(I)+\delta_{2,n}.
\end{equation*}
The complete wild-conductor table is
\begin{center}
\small
\renewcommand{\arraystretch}{1.2}
\begin{tabular}{@{}lll@{}}
\toprule
$I$ & base exponent $\mathfrak f_2$ & $\delta_{2,n}$ \\
\midrule
$C_3$ & $2$ & $0$ \\
$C_2,C_6$ & $4$ & $\epsilon_n(C_2)$ \\
$C_2,C_6$ & $6$ & $2\epsilon_n(C_2)$ \\
$C_4$ & $8$ & $2\epsilon_n(C_4)+\epsilon_n(C_2)$ \\
$Q_8$ & $5$ & $\epsilon_n(Q_8)+\tfrac12\epsilon_n(C_2)$ \\
$Q_8$ & $6$ & $\epsilon_n(Q_8)+\epsilon_n(C_2)$ \\
$Q_8$ & $8$ & $\epsilon_n(Q_8)+\epsilon_n(C_4)+\epsilon_n(C_2)$ \\
$\SL_2(\mathbb F_3)$ & $3$ &
$\tfrac13\epsilon_n(Q_8)+\tfrac16\epsilon_n(C_2)$ \\
$\SL_2(\mathbb F_3)$ & $4$ &
$\tfrac13\epsilon_n(Q_8)+\tfrac23\epsilon_n(C_2)$ \\
$\SL_2(\mathbb F_3)$ & $6$ &
$\tfrac13\epsilon_n(Q_8)+\tfrac53\epsilon_n(C_2)$ \\
$\SL_2(\mathbb F_3)$ & $7$ &
$\tfrac53\epsilon_n(Q_8)+\tfrac56\epsilon_n(C_2)$ \\
\bottomrule
\end{tabular}
\end{center}
Every expression in the last column is an integer.
This is \cite[Tables~2 and 3]{MartinWatkins}; those tables give the wild
exponent $\delta_{2,1}$, whereas $\mathfrak f_2$ here is the total conductor
$\epsilon_1(I)+\delta_{2,1}$.

Let $w_2=W_2(E)$.  For $I=C_2,C_4,C_6$ and odd $n$,
\begin{equation}\label{eq:p2-cyclic-root}
 W_{2,n}=
 \begin{cases}
 w_2,&n\equiv1\pmod4,\\
 1,&n\equiv3\pmod4.
 \end{cases}
\end{equation}
For $I=C_3$, the sequence for
$n\equiv1,3,5,7,9,11\pmod{12}$ is
\begin{equation}\label{eq:p2-C3-root}
 w_2,\ -1,\ -w_2,\ -1,\ w_2,\ 1.
\end{equation}
For $I=Q_8$ or $\SL_2(\mathbb F_3)$,
\begin{equation}\label{eq:p2-exotic-root}
 W_{2,n}=
 \begin{cases}
 w_2,&n\equiv1,5\pmod8,\\
 (-1)^{\mathfrak f_2},&n\equiv3\pmod8,\\
 1,&n\equiv7\pmod8.
 \end{cases}
\end{equation}
The $Q_8$ values of $\mathfrak f_2$ are $5,6,8$, and the
$\SL_2(\mathbb F_3)$ values are $3,4,6,7$.  In the latter case
\eqref{eq:p2-exotic-root} is obtained by passing to the non-Galois tame cubic
extension used in \cite{DMW}, reducing to $Q_8$, and then using the conductor
parity in the table above.  The value $w_2$ itself requires the minimal-model
data.
The $Q_8$ formulas are \cite[Propositions~8.1 and 8.2]{DMW}; the
$\SL_2(\mathbb F_3)$ reduction is \cite[Section~10]{DMW}.

\section{CM curves and the archimedean place}\label{app:cm-infinity}

\subsection{CM decomposition}

\begin{proposition}[Symmetric powers in the CM case]\label{prop:CM}
Suppose that $E$ has complex multiplication by an imaginary quadratic field
$K$.  In the unitary normalization write
\begin{equation*}
 \pi_E^u=\AI_{K/\Q}(\mu),
 \qquad
 \mu^c=\mu^{-1},
\end{equation*}
and let $\chi_K$ be the quadratic character of $K/\Q$.  For $r\geq0$ the
odd formula is
\begin{align}
 \Sym^{2r+1}\pi_E^u
 &=
 \boxplus_{j=0}^{r}
 \AI_{K/\Q}\left(\mu^{2r+1-2j}\right),
 \label{eq:CM-odd}\\
 \intertext{and, for $r\geq1$, the even formula is}
 \Sym^{2r}\pi_E^u
 &=
 \left(
 \boxplus_{j=0}^{r-1}
 \AI_{K/\Q}\left(\mu^{2r-2j}\right)
 \right)
 \boxplus\chi_K^r.
 \label{eq:CM-even}
\end{align}
Also $\Sym^0\pi_E^u=1$; thus no empty isobaric sum is required.
\end{proposition}

\begin{proof}
Restrict the induced two-dimensional parameter to $W_K$.  It is
$\mu\oplus\mu^{-1}$.  Pair the weights $\mu^{n-2j}$ and
$\mu^{-(n-2j)}$.  Each nonzero pair induces to the displayed
two-dimensional automorphic induction.  For even $n=2r$, the middle weight is
one-dimensional and descends as $\chi_K^r$.
\end{proof}

The one-dimensional middle term in \eqref{eq:CM-even} must not be written as
an automorphic induction.  For $R_{a,b}$, a zeta factor occurs internally in
the CM decomposition exactly when $4\mid a$.  It is
\begin{equation}\label{eq:CM-zeta}
 \zeta\left(s-b-\frac a2\right).
\end{equation}
When $a\equiv2\pmod4$, the terminal factor is instead
\begin{equation}\label{eq:CM-quadratic}
 L\left(s-b-\frac a2,\chi_K\right).
\end{equation}

\subsection{The archimedean factor}

Define
\begin{equation*}
 \GammaR(s)=\pi^{-s/2}\Gamma(s/2),
 \qquad
 \GammaC(s)=2(2\pi)^{-s}\Gamma(s).
\end{equation*}
The unitary archimedean parameter is
\begin{equation*}
 \phi_{E,\infty}
 =\Ind_{W_\C}^{W_\mathbb R}(z/|z|).
\end{equation*}

\begin{proposition}\label{prop:infinity}
For $n=2r+1$,
\begin{equation}\label{eq:infinity-odd-unitary}
 \Sym^n\phi_{E,\infty}
 =
 \bigoplus_{j=0}^{r}
 \Ind_{W_\C}^{W_\mathbb R}
 \left((z/|z|)^{n-2j}\right),
\end{equation}
and
\begin{equation}\label{eq:infinity-odd-motivic}
 L_\infty(s,\Sym^nE)
 =
 \prod_{j=0}^{r}\GammaC(s-j).
\end{equation}
For $n=2r$,
\begin{equation}\label{eq:infinity-even-unitary}
 \Sym^n\phi_{E,\infty}
 =
 \left(
 \bigoplus_{j=0}^{r-1}
 \Ind_{W_\C}^{W_\mathbb R}
 ((z/|z|)^{2r-2j})
 \right)
 \oplus\operatorname{sgn}^{\,r},
\end{equation}
and
\begin{equation}\label{eq:infinity-even-motivic}
 L_\infty(s,\Sym^nE)
 =
 \left(\prod_{j=0}^{r-1}\GammaC(s-j)\right)
 \GammaR(s-r+(r\bmod2)).
\end{equation}
For $R_{a,b}$ the factor is $L_\infty(s-b,\Sym^aE)$.
\end{proposition}

\begin{proof}
Pair in $\Sym^n\phi_{E,\infty}$ the characters
$(z/|z|)^{n-2j}$ and $(z/|z|)^{-(n-2j)}$.  Each nonzero pair gives the
displayed two-dimensional induction.  When $n=2r$, the middle character
extends to $W_\mathbb R$ as $\operatorname{sgn}^{r}$.  In the unitary
normalization,
\begin{align*}
 L_\infty\left(s,\Ind_{W_\C}^{W_\mathbb R}(z/|z|)^m\right)
 &=\GammaC(s+m/2)\quad(m\geq1),\\
 L_\infty(s,\operatorname{sgn}^{r})
 &=\GammaR(s+(r\bmod2)).
\end{align*}
More precisely,
$\Sym^n\tau_{E,\infty}=\Sym^n\phi_{E,\infty}
\otimes\omega_\infty^{-n/2}$.  Hence replacing $s$ by $s-n/2$ gives
\eqref{eq:infinity-odd-motivic} and
\eqref{eq:infinity-even-motivic}.  These are the standard local Langlands
conventions of \cite{GetzHahn}; in particular, the factor $2$ in
$\GammaC(s)=2(2\pi)^{-s}\Gamma(s)$ is intentional.  The archimedean root
number follows by multiplying the epsilon factors of the same summands,
giving \eqref{eq:infinity-root} below.
\end{proof}

The archimedean root number is $1$ for even $n$.  For $n=2r+1$, it is
\begin{equation}\label{eq:infinity-root}
 W_{\infty,n}=(-1)^{(r+1)(r+2)/2}.
\end{equation}
Equivalently, $R_{a,b}$ has Hodge types
\begin{equation}\label{eq:Hodge-types}
 (a-j+b,j+b),\qquad 0\leq j\leq a.
\end{equation}
It is pure of weight $a+2b$.  A homogeneous expression of degree $D$ has
functional-equation center $(D+1)/2$.  An isobaric sum involving different
degrees is mixed and has no single motivic center, although every constituent
can be unitary normalized separately.

\section*{Acknowledgements}

The author thanks Jayce Getz and Chen Wan for helpful information about
isobaric sums.  The author also thanks Bin Guan, Xinchen Miao, and Ruichen Xu
for useful discussions, and Bingrong Huang and Philippe Michel for their
encouragement.

\end{document}